\documentclass{amsart}

\usepackage{amssymb}

\theoremstyle{plain}

\numberwithin{equation}{section}

\newtheorem{lem}[equation]{Lemma}

\newtheorem*{teo}{Theorem A}
\newtheorem*{teoremab}{Theorem B}
\newtheorem*{teoremac}{Theorem C} 
\newtheorem{hypo}[equation]{Hypothesis}

\newcommand{\Cl}{\operatorname{Cl}}
\newcommand{\Aff}{\operatorname{Aff}}
\newtheorem{prop}[equation]{Proposition} 

\begin{document}	
\title{Conjugacy classes and finite $p$-groups}

\author{Edith Adan-Bante}

\address{University of Southern Mississippi Gulf Coast, 730 East Beach Boulevard,
 Long Beach MS 39560} 

\email{Edith.Bante@usm.edu}

\keywords{Conjugacy class, $p$-groups, products}

\subjclass{20D15}

\date{2004}

\begin{abstract} Let $G$ be a finite $p$-group, where $p$ is a prime number
and $a\in G$. Denote by $\Cl(a)=\{gag^{-1}\mid g\in G\}$ 
the conjugacy class of $a$ in $G$. Assume that $|\Cl(a)|=p^n$.
Then  $\Cl(a)\Cl(a^{-1})=\{xy\mid x\in \Cl(a), y\in \Cl(a^{-1})\}$  
 is the union of at least $n(p-1)+1$ distinct conjugacy 
classes of $G$.
\end{abstract}
 
\maketitle

\begin{section}{Introduction}

Let $G$ be a finite group. Denote by $\Cl(a)=\{gag^{-1}\mid g\in G\}$ 
the conjugacy class of $a$ in $G$, and by $|\Cl(a)|$ the   size of
$\Cl(a)$.
   If the subset  $X$ of $G$ is $G$-invariant, i.e  $X^g=\{x^g\mid x\in X\}=X$ for
all $g\in G$, then $X$ is the 
union of $m$ distinct conjugacy classes of $G$, for some integer $m$.  
Set $\eta(X)=m$.

Given any conjugacy classes $\Cl(a)$ and $\Cl(b)$, we can check that 
 the product  $\Cl(a)\Cl(b)=\{xy\mid x\in \Cl(a), y \in \Cl(b)\}$ is a 
 $G$-invariant set.
In this note, we will explore the relation between $|\Cl(a)|$ and 
$\eta(\Cl(a)\Cl(a^{-1}))$. Those results are the equivalent in conjugacy classes
as some of the ones in irreducible characters in \cite{edith} and \cite{edith2}.
 
In Theorem A of \cite{edith2}, it is proved that if $G$ is a $p$-group,  
$\chi$ is an irreducible character with degree $p^n$, then the product
$\chi\overline{\chi}$ of $\chi$ with its complex conjugate $\overline{\chi}$
 is the sum of
at least $2n(p-1)+1$ distinct irreducible characters. The following is
the equivalent for conjugacy classes  

\begin{teo} Let $G$ be a finite $p$-group and $a\in G$. Assume that 
 $|\Cl(a)|=p^n$. Then the
product $\Cl(a)\Cl(a^{-1})$ of the conjugacy class of $a$ in $G$ and 
the conjugacy class of the inverse of $a$ in $G$, is the union 
of at least $n(p-1)+1$ distinct conjugacy classes of $G$, 
i.e. $\eta(\Cl(a)\Cl(a^{-1}))\geq n(p-1)+1$.
\end{teo}

In Proposition \ref{examplenp}, it is shown 
 that for every prime $p$ and every 
integer $n\geq 0$, 
there exist a $p$-group $G$ and a conjugacy class $\Cl(a)$ of $G$ such that
$|\Cl(a)|=p^n$ and $\eta(\Cl(a)\Cl(a^{-1}))=n(p-1)+1$. Thus the bound in Theorem A is 
optimal. 

An application of Theorem A is the following
\begin{teoremab}
Let $n$ be a positive integer. Then there exists a finite set $S_n$ of positive integers
such that for any nilpotent group $G$ and any conjugacy class $\Cl(a)$ of $G$ with
$\eta(\Cl(a)\Cl(a^{-1}))\leq n$,  we have that 
$$|\Cl(a)|\in S_n.$$
\end{teoremab}

In Proposition \ref{notforsuper}, we prove that given any prime $p$, there exist
a supersolvable group and a conjugacy class $\Cl(a)$ of G with $|\Cl(a)|=p$
and $\eta(\Cl(a)\Cl(a^{-1}))=2$. Thus the previous result does not
remain true assuming the weaker hypothesis that the groups are supersovable.
Theorem A is the equivalent in conjugacy classes of Theorem B of \cite{edith}.

\begin{teoremac}
Let $p$ be a prime number. Let $G$ be a finite $p$-group and $\Cl(a)$ be a conjugacy 
class of $G$. Then one of the following holds:

i) $|\Cl(a)|=1$ and $\eta(\Cl(a)\Cl(a^{-1}))=1$.

ii) $|\Cl(a)|=p$ and $\eta(\Cl(a)\Cl(a^{-1}))=p$.

iii) $|\Cl(a)|\geq p^2$ and $\eta(\Cl(a)\Cl(a^{-1}))\geq 2p-1$.
\end{teoremac}

Given a fix prime $p>2$, observe that Theorem C implies that there are ``gaps" 
 among the possible 
values that $\eta(\Cl(a)\Cl(a^{-1}))$ can take for any finite 
$p$-group and any 
conjugacy class $\Cl(a)$ in $G$.
 The previous result is the equivalent in 
conjugacy classes of Theorem B of \cite{edith2}.

{\bf Acknowledgment.} I would like to thank 
Manoj Kumar for bringing to my attention 
products of conjugacy classes. 
I also want to thank Professor Everett C. Dade for useful advise and 
corrections.

\end{section}

\begin{section}{Proof of Theorem A}
{\bf Notation}. Let $G$ be a finite $p$-group and $N$ be a normal subgroup of $G$.
Denote by $ \overline{a}$ the element in $G/N$ that contains $a$. Thus 
$\Cl( \overline{a})$ is the conjugacy class of 
$ \overline{a}$ in $G/N$. 
 
\begin{lem}\label{lemma1}
Let $G$ be a finite $p$-group and $N$ be a normal subgroup of $G$. 
Let $a$ and $b$ be elements
of $G$. Then

i) $\Cl( \overline{a})\Cl( \overline{b})$ is a $G$-invariant set.
If $\Cl( \overline{a})\cap \Cl( \overline{b})= \emptyset$ then $\Cl(a)\cap \Cl(b)= \emptyset$.
 Thus $\eta(\Cl( \overline{a})\Cl( \overline{b}))\leq  \eta(\Cl(a)\Cl(b))$.
 
ii) If, in addition,  $|N|=p$, then either $|\Cl(\overline{a})|=|\Cl(a)|$ or 
 $|\Cl(\overline{a})|=\frac{|\Cl(a)|}{p}$. Furthermore, if $|\Cl(\overline{a})|=\frac{|\Cl(a)|}{p}$, 
 then $\eta(\Cl(a)\Cl(a^{-1}))\geq \eta(\Cl(\overline{a})\Cl(\overline{a}^{-1}))+(p-1)$.
 
 \end{lem}
\begin{proof}
{\bf i)}
Clearly if $a=gbg^{-1}$, then $\overline{a}=\overline{g}\overline{b}(\overline{g})^{-1}$.
Thus  if $\Cl(\overline{a})\cap \Cl(\overline{b})= \emptyset$ then $\Cl(a)\cap \Cl(b)= \emptyset$.
Therefore $\eta(\Cl(\overline{a})\Cl(\overline{b}))\leq  \eta(\Cl(a)\Cl(b))$.

{\bf ii)} 
 Since $N$ is normal, $|N|=p$ and  $G$ is a $p$-group, then $N$ is contained in 
 the center ${\bf Z}(G)$ of $G$. Thus given any $n\in N$, $\Cl(n)=\{n\}$. 
 
Suppose that $|\Cl(\overline{a})|\neq |\Cl(a)|$. Since $|N|=p$, we have that
$|\Cl(a)|\leq p |\Cl( \overline{a})|$. Therefore 
$\frac{|\Cl(a)|}{p}\leq |\Cl(\overline{a})|\leq |\Cl(a)|$. Thus $|\Cl(\overline{a})|=\frac{|\Cl(a)|}{p}$
since $G$ is a $p$-group   and $|\Cl(\overline{a})|$ divides $|G/N|$.

If $|\Cl(\overline{a})|=\frac{|\Cl(a)|}{p}$, then given any $x\in \Cl(a)$ and any 
$n\in N$, we have that $nx\in \Cl(a)$. Thus $n=nx(x^{-1})\in \Cl(a)\Cl(a^{-1})$
for any $n\in N$. Therefore
 $N\leq \Cl(a)\Cl(a^{-1})$ and ii) follows. 
\end{proof} 
 \begin{proof}[Proof of Theorem A]
We are going to use  induction on the order of $G$. Let
 $N$ be a normal subgroup of $G$
of
order $p$. Observe such group exists since $G$ is a $p$-group. 
Let $|\Cl(\overline{a})|=p^m$.  
Since $|G/N|<|G|$, 
by induction we have that 
 $\eta(\Cl(\overline{a})\Cl(\overline{a}^{-1}))\geq m(p-1)+1$. 
If $|\Cl(\overline{a})|=|\Cl(a)|$, i.e if $m=n$,
 then by Lemma \ref{lemma1} i) we have that
 $$\eta(\Cl(a)\Cl(a^{-1}))\geq \eta(\Cl(\overline{a})\Cl(\overline{a}^{-1}))\geq m(p-1)+1=n(p-1)+1. $$
We may assume then that  $|\Cl(\overline{a})|\neq|\Cl(a)|$.
By Lemma \ref{lemma1} ii), we have that  $m=n-1$ and 
\begin{equation*}
\begin{split}
\eta(\Cl(a)\Cl(a^{-1})) &\geq  \eta(\Cl(\overline{a})\Cl(\overline{a}^{-1}))+(p-1)\\
                       &=  (n-1)(p-1)+ 1+ (p-1)=n(p-1)+1.                       
\end{split}
\end{equation*}
\end{proof}
 \end{section} 

\begin{section}{Proof of Theorem B}
\begin{proof}[Proof of Theorem B]
Let 
\begin{equation*}
S_n=\{\Pi(p_i)^{t_i}\mid  p_i \mbox{ is a prime number for all } i, \, t_i\geq 0 
\mbox{ and } t_i(p_i-1)+1\leq n \}.
\end{equation*}
Observe that the set $S_n$ is a finite set of positive integers since $0\leq t_i\leq n $ 
and if $t_i>0$ then $p_i\leq n$.

Let $\{p_1, \ldots, p_r\}$ be the set of distinct prime divisors of $|G|$.
 For $i=1, \ldots,r$, 
let $P_i$ be the Sylow $p_i$-subgroup of $G$. Observe that $a=\Pi_{i=1}^r a_i$, for some  
$a_i\in P_i$ for $i=1,\ldots, r$. Since $G$ is nilpotent, we have 
that  $\Cl(a)=\Pi_{i=1}^r \Cl(a_i)$, where 
$\Cl(a_i)$ is the conjugacy class of $a_i$ in $P_i$, for $i=1,\ldots,r$.
  Let $m_i= \eta (\Cl(a_i)\Cl(a_i^{-1}))$. 
Observe that $m_i\leq n$ and
$$|\Cl(a)|=\Pi_{i=1}^r |\Cl(a_i)|.$$
We can check that $$\eta(\Cl(a)\Cl(a^{-1}))=\Pi_{i=1}^r \eta(\Cl(a_i)\Cl(a_i^{-1})).$$
For each $i$, let $|\Cl(a_i)|=p_i^{t_i}$. 
Since $\Cl(a_i)$ is the conjugacy class of $a_i$ in the $p_i$-group $P_i$, 
 by Theorem A we have that 
$m_i=\eta(\Cl(a_i)\Cl(a_i^{-1}))\geq t_i(p_i-1)+1$. Thus $n\geq t_i(p_i-1)+1$.
Therefore 
$|\Cl(a)|=\Pi_{i=1}^r |\Cl(a_i)|= \Pi_{i=1}^r p_i^{t_i}\in S_n.$ 
\end{proof} 
\end{section}
\begin{section}{Proof of Theorem C}
\begin{lem}\label{basicp}
Let $G$ be a finite $p$-group and $\Cl(a)$ be a conjugacy class of $G$ 
with $|\Cl(a)|=p$.
Then one of the following holds:

i) $\Cl(a)=\{az|z \in Z\}$ for some  subgroup $Z$ of the center ${\bf Z}(G)$
of $G$. Therefore
$\Cl(a)\Cl(a^{-1})=Z$ and $\eta(\Cl(a)\Cl(a^{-1}))=p$.

ii) $\Cl(a)\Cl(a^{-1})$ is the union of $p-1$ distinct conjugacy classes 
of size $p$ and 
the class $\Cl(e)=\{e\}$. Therefore $\eta(\Cl(a)\Cl(a^{-1}))=p$.
\end{lem}
\begin{proof}
Observe that  if $z\in \Cl(a)\Cl(a^{-1})$ and $|\Cl(z)|=1$, then 
$z$ is in the center  $\textbf{Z}(G)$ of $G$. Since $z=a^g a^{-1}$ for some $g\in G$ and 
$z\in {\bf Z}(G)$, $z^{i}\in \Cl(a)\Cl(a^{-1})$
and $a^{g^i}=a z^i$  for all integer $i$. Thus $<z>\leq \Cl(a)\Cl(a^{-1})$.
Set $Z=<z>$.

\textbf{i)}
If $z\neq e$, it follows that $|Z|\geq p$. Since $|\Cl(a)|=p$ and $a^{g^i}=a z^i$  for all integer $i$, we have that $\Cl(a)=\{az|z \in Z\}$ and $|Z|=p$.
 Since $Z$ is contained in
${\bf Z}(G)$, then $\Cl(a)\Cl(a^{-1})=Z$ and $\eta(\Cl(a)\Cl(a^{-1}))=p$.

\textbf{ii)}
We may assume now that if  $z\in \Cl(a)\Cl(a^{-1})$ and $\Cl(z)=1$, then 
$z=e$. Thus all the conjugacy classes different from $\Cl(e)$ are of size $p$.
Observe that $a^g (a^{-1})^g =e$ for all $g \in G$.
Thus $|\Cl(a)\Cl(a^{-1})|\leq p^2- p+1 =(p-1)p+1$.
Therefore by Theorem A it follows that $\Cl(a)\Cl(a^{-1})$ is the 
union of  $p-1$ distinct conjugacy classes
of size $p$ and $\Cl(e)$. 
\end{proof}

\textbf{Remark.} Let $p$ be a prime number.

{\bf a)}
Let $G$ be an extra special group of order $p^3$ and exponent $p$.
We can check that given any $a\in G$, where $a$ is not in the center of $G$,
then $\Cl(a)\Cl(a^{-1})={\bf Z}(G)$ and thus Lemma \ref{basicp} i) occurs.

{\bf b)}  
 Let $G$ be the wreath product of a cyclic group $C_{p^2}$ of order $p^2$ 
 by a cyclic group $C_p$
of order $p$. Thus $|G|= p^{2p+1}$. 
Let $a=(c,e,\ldots, e)$ in $G$, where $c\in C_{p^2}$ has order $p^2$.
Observe that $a^{-1}\neq a$.  
Observe also that  $$\Cl(a)=\{(c,e,\ldots, e), (e,c,\ldots, e), \ldots, (e,e,\ldots, e,c)\}.$$
Thus $|\Cl(a)|=p$. Let  
 $b_i=(c, e, \ldots , c^{-1}, e, \ldots, e)$, where $c^{-1}$ is in the 
$i$-position for $i=0, \ldots, p-1$, i.e  $b_0=(c c^{-1},\ldots, e)=(e,e,\ldots, e)$, $b_1=(c,c^{-1}, e,\ldots,e)$ and so for. Observe that $\Cl(b_0)=\Cl((e,e,\ldots, e))$
has class size 1. We can check that $\Cl(b_i)$ has size $p$ for $i=1,\ldots, p-1$.
  Since $c\neq c^{-1}$, then 
$\Cl(b_i)\cap \Cl(b_j)=\emptyset$ if $i\neq j$ and $i, j=0, \ldots, p-1$. 
 Observe that $\Cl(a)\Cl(a^{-1})= \cup_{i=0}^{p-1} \Cl(b_i)$.
 Thus $\Cl(a)\Cl(a^{-1})$
  is the 
union of  a conjugacy class of size 1, namely $\Cl(b_0)$
and $p-1$ distinct conjugacy classes
of size $p$, namely $\Cl(b_i)$ for $i=1,\ldots p-1$. We 
conclude that given any prime $p$, there
exist some group $G$ and some conjugacy class $\Cl(a)$  of $\,G$ satisfying
the condition in case ii) of Lemma \ref{basicp}. 

\begin{proof}[Proof of Theorem C]
If $|\Cl(a)|=\{a\}$, then $\Cl(a)\Cl(a^{-1})=\{e\}$ and so 
i) holds. Lemma \ref{basicp} implies ii) and iii) follows 
from Theorem A.
\end{proof}
\end{section}

\begin{section}{Examples}
  \begin{lem}\label{construction}
 Let $G_0$ be a $p$-group  and $\Cl(g_0)$ be the conjugacy class containing $g_0\in G_0$.
 Assume that $\Cl(g_0)\neq \Cl(g_0^{-1})$.
Let $N=G_0\times G_0\times \cdots\times G_0$ be the direct product of $p$-copies of $\,G_0$.
  Let
$C=<c>$ be a cyclic group of order $p$. 
 Observe that $C$ acts on $N$ by
 \begin{equation}\label{actionc}
 c:(n_0, n_1,\ldots,n_{p-1})\mapsto(n_{p-1},n_0, \dots, n_{p-2})
 \end{equation}
 \noindent for any $(n_0, n_1, \ldots,n_{p-1})\in N$.
 
Let $G$ be the semidirect product of $N$ and $C$, i.e $G$ is the wreath product of $G_0$ and 
$C$.  Set $a=(g_0, e, \ldots, e)$ in $N$, where $e$ is the identity of $G_0$. Then $|\Cl(a)|=p|\Cl(g_0)|$, $\Cl(a)\neq \Cl(a^{-1})$ and 
$\eta(\Cl(a) \Cl(a^{-1}))= \eta(\Cl(g_0) \Cl({g_0}^{-1}))+(p-1)$. 
 \end{lem}
 \begin{proof} 
 Observe that $\Cl(a)\neq \Cl(a^{-1})$ since $\Cl(g_0)\neq \Cl(g_0^{-1})$. 
 
Let $\Cl(g_0) \Cl({g_0}^{-1})=C_1\cup C_2 \cdots \cup C_m$, where $C_1, \ldots, C_m$ 
 are distinct conjugacy classes of $G_0$. Thus $m=\eta(\Cl(g_0)\Cl(g_0^{-1}))$. 
 We can check that  
 the distinct conjugacy classes of  $\Cl(a) \Cl(a^{-1})$ 
 are of the following two types: 
 
 i) $\{(x, e,\ldots,e)^c \mid x\in C_i, c\in C \} $ for
 $i=1, \ldots, m$.
 
 ii)$ \{(x,y,\ldots,  e,e)^c\mid x\in \Cl(g_0), y \in \Cl({g_0}^{-1}), c\in C \}$, 
    $ \{(x,e,y, \ldots,e)^c\mid x\in \Cl(g_0), y \in \Cl({g_0}^{-1}), c\in C \}$ and
    $\{(x,e,\ldots, e,y)^c\mid x\in \Cl(g_0), y \in \Cl({g_0}^{-1}), c\in C \}$.  

Observe that there are $\eta(\Cl(g_0) \Cl({g_0}^{-1}))$ distinct
conjugacy classes of  type i) and 
exactly  $p-1$ distinct conjugacy classes of type ii).  
Thus $\eta(\Cl(a) \Cl(a^{-1}))= \eta(\Cl(g_0) \Cl({g_0}^{-1}))+(p-1)$.
 \end{proof}
 
 \begin{prop}\label{examplenp}
 Given any prime $p$, and any integer $n\geq 0$, there exist a finite
 $p$-group 
 $\,G$ and a conjugacy class $\,\Cl(a)$ of $\,G$ with $|\Cl(a)|=p^n$, $\Cl(a)\neq \Cl(a^{-1})$
  and
 $\eta(\Cl(a) \Cl(a^{-1}))=n(p-1)+1$.
 
\end{prop}
\begin{proof}
Observe that if $G$ is an abelian  group and $a\in G$ has order $p^2$, then $|\Cl(a)|=1$,
$\Cl(a)\neq \Cl(a^{-1})$
and $\eta(\Cl(a) \Cl(a^{-1}))= 1= 0(p-1)+1$. Thus the statement is true for 
$n=0$. Assume by induction that the statement is true for 
$n-1$, i.e. there exist a finite $p$-group $G_0$
and a conjugacy class $\Cl(g_0)$ of $G_0$ with $|\Cl(g_0)|=p^{n-1}$, $\Cl(g_0)\neq \Cl(g_0^{-1})$ and 
 $\eta(\Cl(g_0) \Cl({g_0}^{-1}))=(n-1)(p-1)+1$. Using the notation of 
 Lemma \ref{construction}, we have that 
 \begin{equation*}
\begin{split}
\eta(\Cl(a) \Cl(a^{-1})) &=
 \eta(\Cl(g_0) \Cl({g_0}^{-1}))+(p-1)\\
 &=(n-1)(p-1)+1+ (p-1)= n(p-1)+1.
 \end{split}
\end{equation*}
\noindent Since $|\Cl(a)|=p|\Cl(g_0)|= p\times p^{n-1}=p^n$, the proof is complete.
\end{proof}

\begin{hypo}\label{hypot}	Fix a prime $p$ and let
$F=\{0,1,\ldots,p-1\}$ be the finite field with $p$ elements. Observe that $F$ is also
a vector
space of dimension 1 over itself.
Let 
$ A=\Aff(F)$ be the affine group of $F$.  
Observe that the group $A$ is a cyclic by cyclic group and thus
it is supersolvable. 

Let $C$ be a cyclic group of order $p$. 
 Set $X=F$ and  $K = C^X = \{ f: X \to C \}$. Observe that
$K$ 
is a group via pointwise multiplication, and clearly
$A$ acts on this group (via its action on $X$).
 
 	Let $G$ be the wreath product of $C$ and $A$ relative to $X$,
i.e. $G = K \rtimes A$. We can check  that $G$ is a supersolvable group.
\end{hypo}
\begin{prop}\label{notforsuper}
Assume Hypotheses \ref{hypot}. 
Set $a=(c,e,e,\ldots, ,e)$ in $K$. 
Then $a\in G$, the conjugacy class $\Cl(a)$ of $\,G$  has size $p$ and 
$\Cl(a)\Cl(a^{-1})=\Cl((e,e,\ldots, e))\cup \Cl((c,c^{-1},e,  \ldots, e))$. 
Thus  $\eta(\Cl(a)\Cl(a^{-1}))=2$.

Therefore, given any prime $p$, there exist a supersolvable group $G$  and
a conjugacy class $\Cl(a)$ of $G$ with $|\Cl(a)|=p$ and $\eta(\Cl(a)\Cl(a^{-1}))=2$.
\end{prop}
\begin{proof}
Observe that 
$\Cl(a)=\{(c, e, \ldots,e), (e,c,e,\ldots, e),\ldots, (e,e,\ldots,e,c)\}$. 
Thus $\Cl(a)$ has $p$-elements.
Observe that $$\{(c, c^{-1}, \ldots , e)^y \mid y \in F\setminus\{0\}\}=
 \{(c, c^{-1},e,  \ldots, e), (c,e, c^{-1},   \ldots, e), 
\ldots, (c, e,  \ldots, e, c^{-1})\}.$$ Observe also that 
$$\{(c, c^{-1},e,  \ldots, e)^x\mid x\in F\}=  
 \{(c, c^{-1},e,  \ldots, e),   
\ldots, (e, e,  \ldots, c, c^{-1}), (c^{-1},e, \ldots,e, c)\}.$$
 Thus  
\begin{equation*}
\begin{split}
\Cl((c,c^{-1}, e,\ldots,e))= & \{(c, c^{-1},e,  \ldots, e)^x, 
(c,e, c^{-1},   \ldots, e)^x,\ldots,  \\
& (c,e, \ldots, c^{-1}, e), (c, e,  \ldots, e, c^{-1})^x\mid x\in F \} 
\end{split}
\end{equation*}
\noindent has $(p-1)p=p^2-p$ elements. Since $a^g (a^{-1})^g=(e,\ldots ,e)$,
then $\Cl(a)\Cl(a^{-1})$ has at most $p^2-p+1$ elements. We conclude that
 $$\Cl(a)\Cl(a^{-1})=\Cl((e,e, \ldots ,e))\cup \Cl((c, c^{-1},e, \ldots, e)).$$  
\end{proof}
\end{section}


\begin{thebibliography}{9}

\bibitem{edith} 
E. Adan-Bante, Products of characters and finite $p$-groups, J. Algebra 277 (1), 236-255.
\bibitem{edith2} 
E. Adan-Bante, Products of characters and finite $p$-groups II,  Arch. Math. (82) (2004),
 289-297. 
  
\end{thebibliography}
\end{document}